\documentclass[12 points]{article}
\usepackage{amssymb}
\usepackage{eucal}
\usepackage{amsmath}
\usepackage{amsthm}
\usepackage{graphicx}
\usepackage{float}

\newcommand{\R}{{\mathbb  R}}  \numberwithin{equation}{section} \newtheorem{thm}{\bf
Theorem}[section]
  
  \theoremstyle{remark}
\newtheorem{rem}{\bf
Remark}[section]

\begin{document}

\title{\large\bf Hessian Operators on Constraint Manifolds} \author{Petre Birtea and Dan Com\u{a}nescu \\ {\small Department of Mathematics, West University of Timi\c soara}\\ {\small Bd. V.
P\^ arvan,
No 4, 300223 Timi\c soara, Rom\^ania}\\ {\small birtea@math.uvt.ro, comanescu@math.uvt.ro}\\ }
\date{ } \maketitle

\begin{abstract}

On a constraint manifold  we give an explicit formula for the Hessian matrix of a cost function that involves the Hessian matrix of a prolonged function and the Hessian matrices of the constraint functions. We give an explicit formula for the case of the orthogonal group ${\bf O}(n)$ by using only Euclidean coordinates on $\R^{n^2}$. An optimization problem on ${\bf SO}(3)$ is completely carried out. Its applications to nonlinear stability problems are also analyzed.

\end{abstract}



\section{Introduction}

The nature of critical points of a smooth cost function $G_S:S\rightarrow \R$, where $(S,\tau)$ is a smooth Riemannian manifold, can be  often determined by analyzing the Hessian matrix of the function $G_S$ at the critical points.
In order to compute the Riemannian Hessian one needs to have a good knowledge of the Riemannian geometry of the manifold $S$ such as the affine connection associated with the metric and geodesic lines, see \cite{absil-mahony-sepulchre-1} and \cite{gallot}. Usually, to carry out explicit computations one has to introduce local coordinate systems on the manifold. These elements are often difficult to construct and manipulate on specific examples. A large class of examples have been studied in \cite{absil-mahony-sepulchre-1}, \cite{edelman}, \cite{ferreira} in connection with optimization algorithms. 

A method to bypass the computational difficulties associated with a Riemannian manifold is to embed it in a larger space, usually an Euclidean space, and transfer the computations  into this more simpler space.

In this paper we give a formula for the Hessian matrix of the function $G_S$ that involves the Hessian matrix of an extended function $G$ and the Hessian matrices of the constraint functions. More precisely, let $G_S:S\rightarrow \R$ be a smooth map defined on the manifold $S$. Suppose that  $S$ is a submanifold of a smooth manifold $M$ that is also the preimage of a regular value for a smooth function $\mathbf{F}:=(F_1,\ldots,F_k):M\rightarrow \mathbb{R}^k$, i.e. $S=\mathbf{F}^{-1}(c)$, where $c$ is a regular value of $\mathbf{F}$. Let $G:M\rightarrow \R$ be a prolongation of the function $G_S$. In the case when $(M,g)$ is a finite dimensional Riemannian manifold we can endow the submanifold $S$ with a Riemannian metric $\tau_c$, constructed in \cite{birtea-comanescu}, that is conformal with the induced metric on $S$ by the ambient metric $g$. The gradient of the restricted function $G_S$ with respect to the Riemannian metric $\tau_c$  can be computed using the gradients with respect to the ambient Riemannian metric $g$ of the prolongation function $G$ and the constraint functions $F_1,...,F_k$.

In order to compute the Hessian operator of the cost function $G_S$ we need to take covariant derivatives of the gradient vector field 
$\text{\bf grad}_{g^S_{ind}}G_S$. The covariant derivative on the submanifold $S$ is related to the covariant derivative of the ambient space in the following way: take the covariant derivative in the ambient space of a prolongation of the vector field $\text{\bf grad}_{g^S_{ind}}G_S$ and project this vector field on the tangent space of the submanifold $S$. The new problem is to find a vector field defined on the ambient space that prolongs $\text{\bf grad}_{g^S_{ind}}G_S$. The solution to this problem is given by the standard control vector field introduced in \cite{birtea-comanescu}. 

In Section 3 we apply the formula found in Section 2 to cost functions defined on the orthogonal group ${\bf O}(n)$. In Section 4 we specialize the formula from the previous section to the 2-power cost function considered in \cite{birtea-comanescu-popa}. 

In the last section we relate our main result to a stability problem for equilibrium points of a dynamical system. We give a new interpretation of the stability results using the augmented function technique introduced in \cite{maddocks-1991}, see also \cite{beck-hall} and \cite{wang-xu}. 

In paper \cite{absil-mahony-trumpf} a formula for the Riemannian Hessian has been proved using the orthogonal projection and the Weingarten map for a general submanifold embedded in an Euclidean space. Another construction of a Hessian operator using orthogonal coordinates on the tangent planes of a submanifold embedded in an Euclidean space has been presented in \cite{donoho}.

\section{Construction of the Hessian operator on constraint manifolds}

As discussed in the Introduction, we work in the following setting.
Let $G_S:S\rightarrow \R$ be a smooth map defined on a manifold $S$. Suppose that  $S$ is a submanifold of a smooth manifold $M$ that is also the preimage of a regular value for a smooth function $\mathbf{F}:=(F_1,\ldots,F_k):M\rightarrow \mathbb{R}^k$, i.e. $S=\mathbf{F}^{-1}(c)$, where $c$ is a regular value of $\mathbf{F}$. Let $G:M\rightarrow \R$ be a prolongation of the function $G_S$.  

We recall the construction and the geometry of the standard control vector field introduced in \cite{birtea-comanescu}.
The $r\times s$ Gramian matrix generated by the smooth functions $f_1,...,f_r,g_1,...,g_s:(M,g)\rightarrow \mathbb{R}$ is defined by the formula
\begin{equation}\label{sigma}
\Sigma_{(g_1,...,g_s)}^{(f_1,...,f_r)}=\left[%
\begin{array}{cccc}
  <\text{\bf grad } g_1,\text{\bf grad } f_{1}> & ... & <\text{\bf grad } g_s,\text{\bf grad } f_{1}> \\
  ... & ... & ... \\

  <\text{\bf grad } g_1,\text{\bf grad } f_r> & ... & <\text{\bf grad } g_s,\text{\bf grad } f_r> \\
\end{array}%
\right].
\end{equation}
{\bf The standard control vector field} has the formula
\begin{align}\label{v0}
    \mathbf{v_0} & =\sum_{i=1}^k(-1)^{i+k+1}\det \Sigma_{(F_1,\ldots ,\widehat{F_i},\ldots ,F_k,G)}^{(F_1,\ldots
    ,F_k)}\text{\bf grad }    F_i+\det\Sigma_{(F_1,\ldots ,F_k)}^{(F_1,\ldots ,F_k)}\text{\bf grad } G \nonumber\\
& = \det\Sigma_{(F_1,\ldots ,F_k)}^{(F_1,\ldots ,F_k)}\text{\bf grad } G - \sum_{i=1}^k\det \Sigma_{(F_1,\ldots ,F_{i-1},G, F_{i+1},\dots,F_k)}^{(F_1,\ldots  ,F_{i-1},F_i, F_{i+1},...  ,F_k)}\text{\bf grad }    F_i,
\end{align}
where  $\widehat{\cdot}$ represents the missing term. The vector field ${\bf v}_0$ is tangent to the submanifold $S$ and consequently, the restriction defines a vector field on $S$, i.e. ${{\bf v}_0}_{|S}\in \mathcal{X}(S)$. 
On the submanifold $S$ we can define a Riemannian metric $\tau_c$ such that, see \cite{birtea-comanescu},
$${{\bf v}_0}_{|S}=\text{\bf grad} _{\tau_c}G_S.$$
We make the notation $\Sigma:=\det\Sigma_{(F_1,\ldots ,F_k)}^{(F_1,\ldots ,F_k)}$.  On the submanifold $S$ the restricted function $\Sigma_{|S}$ is everywhere different from zero as the submanifold $S$ is the preimage of a regular value. Also, it has been proved in \cite{birtea-comanescu} that the relation between the Riemannian metric $\tau_c$ defined on $S$ and the induced Riemannian metric $g^S_{ind}$ from the ambient space $(M,g)$ is given by
$$\tau_c=\frac{1}{\Sigma_{|S}}g^S_{ind}.$$

To determine the relation between the gradient of the function $G_S$ with respect to the metric $\tau_c$ and the gradient of the function $G_S$ with respect to the induced metric $g^S_{ind}$ we have the following computation:
\begin{align*}
g^S_{ind}(\text{\bf grad}_{g^S_{ind}}G_S,{\bf w}) & =dG_S({\bf w})=\tau_c(\text{\bf grad}_{\tau_c}G_S,{\bf w}) \\
 & =\frac{1}{\Sigma_{|S}}g^S_{ind}(\text{\bf grad}_{\tau_c}G_S,{\bf w}) \\
& =g^S_{ind}(\frac{1}{\Sigma_{|S}}\text{\bf grad}_{\tau_c}G_S,{\bf w}) ,\,\,\forall {\bf w}\in \mathcal{X}(S).
\end{align*}
Consequently, 
$$\text{\bf grad}_{g^S_{ind}}G_S=\frac{1}{\Sigma_{|S}}{{\bf v}_0}_{|S}.$$
The above equality implies that a prolongation of the vector field $\text{\bf grad}_{g^S_{ind}}G_S$ to the open subset $\Omega=\{ x\in M\,|\,\Sigma(x)\neq 0\}$ of the ambient space $M$ is given by the vector field 
\begin{equation}\label{v0-grad}
\frac{1}{\Sigma}{{\bf v}_0}=\text{\bf grad }G-\sum_{i=1}^k\sigma_i\text{\bf grad }F_i,
\end{equation}
where 
$\sigma_i:\Omega\rightarrow \R$ are defined by
\begin{equation}\label{sigma}
\sigma_i(x):=\frac{\det \Sigma_{(F_1,\ldots ,F_{i-1},G, F_{i+1},\dots,F_k)}^{(F_1,\ldots , F_{i-1},F_i, F_{i+1},...
    ,F_k)}(x)}{\Sigma(x)}.
\end{equation}

{\bf     If $x_0\in S$ is a critical point of the function $G_S$, then the numbers $\sigma_i(x_0)$ are the Lagrange multipliers of the extended function $G$ constraint to submanifold $S$.}
More precisely, a critical point $x_0$ of the constraint function $G_S=G_{|S}$ is an equilibrium point for the standard control vector field ${\bf v}_0$ which implies the equality that gives the Lagrange multipliers
$$\text{\bf grad }G(x_0)=\sum_{i=1}^k\sigma_i(x_0)\text{\bf grad }F_i(x_0).$$
The above Lagrange multipliers are uniquely determined due to regular value condition which implies that ${\bf grad }F_1(x_0),...,{\bf grad }F_k(x_0)$ are linearly independent vectors in $T_{x_0}M$.
\medskip

In what follows we show how the Hessian operator associated to the cost function $G_S:S\rightarrow \R$, where $S$ is endowed with the induced metric $g^S_{ind}$, is related with the Hessian operator of the extended function $G$ and the Hessian operators of the functions $F_i$ that describe the constraint submanifold $S$. By definition, see \cite{absil-mahony-sepulchre-1}, the Hessian operator $\mathcal{H}^{G_S}(x):T_xS\rightarrow T_xS$ is defined by the equality 
$$ \mathcal{H}^{G_S}(x)\cdot {\boldsymbol\eta}_{_x}=\nabla^S_{{\boldsymbol \eta}_{_x}}\text{\bf grad}_{g^S_{ind}}G_S,$$
where $\nabla^S$ is the covariant derivative on the Riemannian manifold $(S,g^S_{ind})$ and ${\boldsymbol \eta}_{_x}\in T_xS$.
The relation between the covariant derivative $\nabla^S$ and the covariant derivative $\nabla$ associated with the ambient Riemannian manifold $(M,g)$ is given by 
$$\nabla ^S_{{\boldsymbol \eta}_{_x}}{\boldsymbol \xi}^S={\bf P}_{T_xS}\nabla_{\tilde{{\boldsymbol \eta}}_{_x}}{\boldsymbol \xi},$$
where ${\bf P}_{T_xS}:T_xM\rightarrow T_xM$ is the orthogonal projection onto $T_xS$  with respect to the scalar product on the tangent space $T_xM$ induced by the ambient metric $g$, ${\boldsymbol \xi}\in \mathcal{X}(M)$ is a prolongation on the ambient  space of the vector field $\boldsymbol{\xi}^S\in \mathcal{X}(S)$, and $\tilde{{\boldsymbol \eta}}_{_x}\in T_xM$ is the vector ${\boldsymbol \eta}_{_x}\in T_xS$ regarded as a vector in the ambient tangent space $T_xM$.

For $x\in S$, using the prolongation given by \eqref{v0-grad} for the vector field $\text{\bf grad}_{g^S_{ind}}G_S$ we obtain:
\begin{align*}
\mathcal{H}^{G_S}(x)\cdot {\boldsymbol\eta}_{_x} & = {\bf P}_{T_xS}\nabla_{\tilde{{\boldsymbol \eta}}_{_x}}\frac{1}{\Sigma}{{\bf v}_0}={\bf P}_{T_xS}\nabla_{\tilde{{\boldsymbol \eta}}_{_x}}\left(\text{\bf grad }G-\sum_{i=1}^k\sigma_i\text{\bf grad }F_i\right) \\
& ={\bf P}_{T_xS}\nabla_{\tilde{{\boldsymbol \eta}}_{_x}}\text{\bf grad }G-\sum_{i=1}^kd\sigma_i(\tilde{{\boldsymbol \eta}}_{_x}){\bf P}_{T_xS}\text{\bf grad }F_i(x)-\sum_{i=1}^k\sigma_i(x){\bf P}_{T_xS}\nabla_{\tilde{{\boldsymbol \eta}}_{_x}}\text{\bf grad }F_i \\
& ={\bf P}_{T_xS}\nabla_{\tilde{{\boldsymbol \eta}}_{_x}}\text{\bf grad }G-\sum_{i=1}^k\sigma_i(x){\bf P}_{T_xS}\nabla_{\tilde{{\boldsymbol \eta}}_{_x}}\text{\bf grad }F_i \\
& = {\bf P}_{T_xS} \mathcal{H}^{G}(x)\cdot \tilde{{\boldsymbol\eta}}_{_x}-\sum_{i=1}^k\sigma_i(x){\bf P}_{T_xS} \mathcal{H}^{F_i}(x)\cdot \tilde{{\boldsymbol\eta}}_{_x}={\bf P}_{T_xS}\left( \mathcal{H}^{G}(x)-\sum_{i=1}^k\sigma_i(x)\mathcal{H}^{F_i}(x)\right)\cdot \tilde{{\boldsymbol\eta}}_{_x}.
\end{align*}
Consequently, the bilinear form associated to the Hessian operator $\mathcal{H}^{G_S}(x)$ is given by 
\begin{align*}
\text{Hess}\,G_S(x)({\boldsymbol \eta}_{_x},{\boldsymbol \xi}_{_x}) & =<\mathcal{H}^{G_S}(x)\cdot {\boldsymbol \eta}_{_x},{\boldsymbol \xi}_{_x}>_{g^S_{ind}} \\
& =< {\bf P}_{T_xS} \mathcal{H}^{G}(x)\cdot \tilde{{\boldsymbol\eta}}_{_x},\tilde{{\boldsymbol \xi}}_{_x}>_{g}-\sum_{i=1}^k\sigma_i(x)<{\bf P}_{T_xS} \mathcal{H}^{F_i}(x)\cdot \tilde{{\boldsymbol\eta}}_{_x},\tilde{{\boldsymbol \xi}}_{_x}>_{g} \\
& =<\mathcal{H}^{G}(x)\cdot \tilde{{\boldsymbol\eta}}_{_x}, {\bf P}_{T_xS} \tilde{{\boldsymbol \xi}}_{_x}>_{g}-\sum_{i=1}^k\sigma_i(x)< \mathcal{H}^{F_i}(x)\cdot \tilde{{\boldsymbol\eta}}_{_x},{\bf P}_{T_xS}\tilde{{\boldsymbol \xi}}_{_x}>_{g} \\
& =<\mathcal{H}^{G}(x)\cdot \tilde{{\boldsymbol\eta}}_{_x},  \tilde{{\boldsymbol \xi}}_{_x}>_{g}-\sum_{i=1}^k\sigma_i(x)< \mathcal{H}^{F_i}(x)\cdot \tilde{{\boldsymbol\eta}}_{_x},\tilde{{\boldsymbol \xi}}_{_x}>_{g}.
\end{align*}

The above considerations lead us to the main result of the paper.

\begin{thm}\label{calcul-hessiana-restransa}
For any $x\in S$, the symmetric covariant tensor associated with the Hessian operator of the cost function $G_S$ has the following formula: 
\begin{equation}\label{Hessian-constraints}
\left[ \normalfont{\text{Hess}}\,{G_S}(x)\right] =\left[\normalfont{\text{Hess}}\,{G}(x)\right]_{|_{T_xS\times T_xS}}-\sum_{i=1}^k\sigma_i(x)\left[ \normalfont{\text{Hess}}\,{F_i}(x)\right]_{|_{T_xS\times T_xS}}.
\end{equation}
\end{thm}

{\bf The above formula is valid for all points $x\in S$, not just for critical points of the cost function $G_S$.}
Choosing a base for the tangent space $T_xS\subset T_xM$, $\{\tilde{\bf f}_{a}\in T_xM\}_{a= \overline{1,\dim S}}$, the components of the Hessian matrix  $\left[ \text{Hess}\,{G_S}(x)\right]$ are given by the following relation between the components of the Hessian matrix of the prolonged function $G$ and the components of the Hessian matrices of the functions $F_i$ that describe the submanifold $S$:
\begin{align}
\left[ \text{Hess}\,{G_S}(x)\right]_{ab} & =\text{Hess}\,{G}(x)( \tilde{\bf f}_{a},  \tilde{\bf f}_{b})-\sum_{i=1}^k\sigma_i(x) \text{Hess}\,{F_i}(x)( \tilde{\bf f}_{a},\tilde{\bf f}_{b}).
\end{align}

{\bf Note that the base for the tangent space  $T_xS$ is computed in the local coordinates of the ambient manifold $M$ and does not imply the knowledge of a local coordinate system on the submanifold $S$.}
We recall that for a $C^2$-function $G:(M,g)\rightarrow \R$, in a local coordinate system on the ambient manifold $M$ we have the formula
$$\left[\text{Hess}\,G(x)\right]_{uv}=<\frac{\partial }{\partial x^u},\mathcal{H}^{G}(x) \frac{\partial }{\partial x^v}>_g=\frac{\partial^2G}{\partial x^{u}\partial x^{v}}(x)-\Gamma_{uv}^{w}(x)\frac{\partial G}{\partial x^{w}}(x),\,\,u,v,w=\overline{1,\dim M},$$
where $\Gamma_{uv}^{w}$ are the Christoffel's symbols associated to the metric $g$.

\section{Hessian operator for cost functions on ${\bf O}(n)$}

We give an explicit formula for the Hessian operator associated to a $C^2$ cost function $G_{{\bf O}(n)}:{\bf O}(n)\rightarrow \R$. The orthogonal group is defined by:
$${\bf O}(n)=\{X\in \mathcal{M}_{n\times n}(\R)\,|\,XX^T=\mathbb{I}_n=X^TX\}.$$
A general element of ${\bf O}(n)$ is represented by an orthonormal frame $\{{\bf x}_i=(x_{i1},...,x_{in})\}_{i=\overline{1,n}}$ in $\R^n$, namely,
$$X=
\left[%
\begin{array}{c}
 {\bf x}_1 \\
  ...  \\
{\bf x}_n
\end{array}%
\right]
=
\left[%
\begin{array}{cccc}
  x_{11} & ... & x_{1n} \\
  ... & ... & ... \\
 x_{n1} & ... & x_{nn}
\end{array}%
\right].
$$
We identify an orthogonal matrix with a vector in $\R^{n^2}$ using the linear map ${\cal I}:\mathcal{M}_{n\times n}(\R)\rightarrow \R^{n^2}$,
\begin{equation}\label{ident}
X\stackrel{{\cal I}}{\longmapsto} \tilde{\bf x}=({\bf x}_1,...,{\bf x}_n).
\end{equation}
Regarded as a subset of $\R^{n^2}$, the orthogonal group ${\bf O}(n)$ can be seen as the preimage $\mathcal{O}(n)\subset  \R^{n^2}$ of the regular value $(\frac{1}{2},...,\frac{1}{2},0,...,0)\in \R^n\times \R^{\frac{(n-1)n}{2}}$ for the constraint functions:
\begin{align}\label{constraints-1}
F_s (\tilde{\bf x}) & =\frac{1}{2}||{\bf x}_s||^2,\,\,s\in \{1,...,n\}, \\
F_{pq} (\tilde{\bf x}) & = <{\bf x}_p,{\bf x}_q>,\,\,1\leq p<q\leq n, \label{constraints-2}
\end{align}
where $<\cdot,\cdot>$ is the Euclidean product in $\R^n$.
Using the above identification we obtain the cost function $G_{\mathcal{O}(n)}:=G_{{\bf O}(n)}\circ {\cal I}:\mathcal{O}(n)\rightarrow \R$.

Starting with the canonical base ${\bf e}_1,...,{\bf e}_n$  in $\R^n$, we obtain the canonical base $\{\tilde{\bf e}_{ij}\}_{i,j=\overline{1,n}}$ in $\R^{n^2}$ by the identification, 
$${\bf e}_i\otimes {\bf e}_j\stackrel{{\cal I}}{\longmapsto} \tilde{\bf e}_{ij}=({\bf 0},...,{\bf 0},{\bf e}_j,{\bf 0},...,{\bf 0}),$$
where the vector ${\bf e}_j$ is on the $i$-th slot and the $n\times n$ matrix ${\bf e}_i\otimes {\bf e}_j$ has $1$ on the $i$-th row and $j$-th column and the rest $0$. 

By direct computations, we have the following formulas for the gradients of the constraint functions:
\begin{align*}
\text{\bf grad }F_s(\tilde{\bf x}) & =\sum_{i=1}^nx_{si}  \tilde{\bf e}_{si}, \\
\text{\bf grad }F_{pq}(\tilde{\bf x}) & =\sum_{i=1}^n(x_{qi}  \tilde{\bf e}_{pi}+x_{pi}  \tilde{\bf e}_{qi}).
\end{align*}
The $n^2\times n^2$ Hessian matrices of the constraint functions are given by:
\begin{align*}
\left[\text{Hess}\,{F_s}(\tilde{\bf x})\right] & =\sum_{j=1}^n\tilde{\bf e}_{sj}\otimes \tilde{\bf e}_{sj}, \\
\left[\text{Hess}\,{F_{pq}}(\tilde{\bf x})\right] & =\sum_{j=1}^n(\tilde{\bf e}_{pj}\otimes \tilde{\bf e}_{qj}+\tilde{\bf e}_{qj}\otimes \tilde{\bf e}_{pj}).
\end{align*}
For a $C^2$ prolongation $G:\R^{n^2}\rightarrow \R$ of the cost function $G_{\mathcal{O}(n)}$ we have the formula for the Hessian matrix
$$\left[\text{Hess}\,{G}(\tilde{\bf x})\right] =\sum_{a,b,c,d=1}^n\frac{\partial^2 G}{\partial x_{ab}x_{cd}}(\tilde{\bf x})\,\tilde{\bf e}_{ab}\otimes \tilde{\bf e}_{cd}.$$

Identifying a rotation $X\in {\bf O}(n)$ with the corresponding point $\tilde{\bf x}\in \mathcal{O}(n)$ and substituting in \eqref{Hessian-constraints} the formulas \eqref{sigma-uri} (see Annexe), we have the following formula for the restricted Hessian matrix:  
\begin{align}\label{Hessian-O(n)}
\left[ \text{Hess}\,{G_{\mathcal{O}(n)}}(\tilde{\bf x})\right] & =\left[ \text{Hess}\,{G}(\tilde{\bf x})\right]_{|_{T_{\tilde{\bf x}}{\mathcal{O}}(n)\times T_{\tilde{\bf x}}{\mathcal{O}}(n)}}-\sum_{s=1}^n<\frac{\partial G}{\partial {\bf x}_s}(\tilde{\bf x}),{\bf x}_s> \left[ \text{Hess}\,{F_s}(\tilde{\bf x})\right]_{|_{T_{\tilde{\bf x}}{\mathcal{O}}(n)\times T_{\tilde{\bf x}}{\mathcal{O}}(n)}} \nonumber \\
& -\frac{1}{2}\sum_{1\leq p<q\leq n}(<\frac{\partial G}{\partial {\bf x}_p}(\tilde{\bf x}),{\bf x}_q>+<\frac{\partial G}{\partial {\bf x}_q}(\tilde{\bf x}),{\bf x}_p>) \left[ \text{Hess}\,{F_{pq}}(\tilde{\bf x})\right]_{|_{T_{\tilde{\bf x}}{\mathcal{O}}(n)\times T_{\tilde{\bf x}}{\mathcal{O}}(n)}},
\end{align}
where we denote by $\frac{\partial G}{\partial {\bf x}_i}(\tilde{\bf x}):=\sum_{j=1}^n\frac{\partial G}{\partial {x}_{ij}}(\tilde{\bf x}){\bf e}_{j}$, which is a vector in $\R^n$.

Equivalently, 
\begin{align*}
\left[ \text{Hess}\,{G_{\mathcal{O}(n)}}(\tilde{\bf x})\right] & =\left[ \text{Hess}\,{G}(\tilde{\bf x})\right]_{|_{T_{\tilde{\bf x}}{\mathcal{O}}(n)\times T_{\tilde{\bf x}}{\mathcal{O}}(n)}}-\sum_{i,j,s=1}^nx_{si}\frac{\partial G}{\partial x_{si}}(\tilde{\bf x})(\tilde{\bf e}_{sj}\otimes \tilde{\bf e}_{sj})_{|_{T_{\tilde{\bf x}}{\mathcal{O}}(n)\times T_{\tilde{\bf x}}{\mathcal{O}}(n)}} \\
& - \frac{1}{2} \sum_{1\leq p<q\leq n}\sum_{i,j=1}^n(x_{qi}\frac{\partial G}{\partial x_{pi}}(\tilde{\bf x})+x_{pi}\frac{\partial G}{\partial x_{qi}}(\tilde{\bf x}))(\tilde{\bf e}_{pj}\otimes \tilde{\bf e}_{qj}+\tilde{\bf e}_{qj}\otimes \tilde{\bf e}_{pj})_{|_{T_{\tilde{\bf x}}{\mathcal{O}}(n)\times T_{\tilde{\bf x}}{\mathcal{O}}(n)}}.
\end{align*}

To write the above formula in a more explicit form we need to choose a base for the tangent space $$T_X{\bf O}(n)=\{X\Omega\ |\ \Omega=-\Omega^T\}.$$ 
It is equivalent to choose a base for $n\times n$ skew-symmetric  matrices. We consider the following base:
\begin{equation}\label{Omega-100}
\Omega_{\alpha\beta}=(-1)^{\alpha+\beta}({\bf e}_{\alpha}\otimes{\bf e}_{\beta}-{\bf e}_{\beta}\otimes{\bf e}_{\alpha}),\,\,1\leq \alpha<\beta\leq n.
\end{equation}
Consequently, a base for the tangent space $T_X{\bf O}(n)$ is given by (see Annexe \eqref{calcul-baza-tangenta}):
$$X\Omega_{\alpha\beta}= (-1)^{\alpha+\beta}\sum_{i=1}^n(x_{i\alpha}{\bf e}_i\otimes{\bf e}_{\beta}-x_{i\beta}{\bf e}_i\otimes{\bf e}_{\alpha}),\,\,1\leq \alpha<\beta\leq n.$$
As ${\cal I}$ is a linear map its differential $d_{X}{\cal I}$ equals ${\cal I}$, we obtain the following base for $T_{\tilde{\bf x}}{\mathcal O}(n)$:
$$\tilde{\omega}_{\alpha\beta}(\tilde{\bf x})= (-1)^{\alpha+\beta}\sum_{i=1}^n(x_{i\alpha}\tilde{\bf e}_{i\beta}-x_{i\beta}\tilde{\bf e}_{i\alpha}),\,\,1\leq \alpha<\beta\leq n.$$

In the base for $T_{\tilde{\bf x}}{\mathcal O}(n)$ chosen as above we have the following formula of the element $(\gamma\tau)(\alpha\beta)$, $1\leq \gamma<\tau\leq n$ and $1\leq \alpha<\beta\leq n$, of the Hessian matrix for the cost function $G_S$:
 \begin{align*}
& \left[ \text{Hess} \,{G_{\mathcal{O}(n)}} (\tilde{\bf x})\right]_{(\gamma\tau)(\alpha\beta)} =\\ 
& =\epsilon\sum_{a,c=1}^n( x_{a\gamma}x_{c\alpha}\frac{\partial^2 G}{\partial x_{a\tau}x_{c\beta}}(\tilde{\bf x})- x_{a\gamma}x_{c\beta}\frac{\partial^2 G}{\partial x_{a\tau}x_{c\alpha}}(\tilde{\bf x}) - x_{a\tau}x_{c\alpha}\frac{\partial^2 G}{\partial x_{a\gamma}x_{c\beta}}(\tilde{\bf x}) + x_{a\tau}x_{c\beta}\frac{\partial^2 G}{\partial x_{a\gamma}x_{c\alpha}}(\tilde{\bf x}))
  \\
& - \epsilon\sum_{i,s=1}^nx_{si}(x_{s\tau}x_{s\beta}\delta_{\gamma\alpha}-x_{s\gamma}x_{s\beta}\delta_{\tau\alpha}-x_{s\tau}x_{s\alpha}\delta_{\gamma\beta}+x_{s\gamma}x_{s\alpha}\delta_{\tau\beta})\frac{\partial G}{\partial x_{si}}(\tilde{\bf x}) \\
& - \frac{\epsilon}{2} \sum_{1\leq p<q\leq n}\sum_{i=1}^n(x_{qi}\frac{\partial G}{\partial x_{pi}}(\tilde{\bf x})+x_{pi}\frac{\partial G}{\partial x_{qi}}(\tilde{\bf x}))
(x_{p\tau}x_{q\beta}\delta_{\gamma\alpha}-x_{p\gamma}x_{q\beta}\delta_{\tau\alpha}
 -x_{p\tau}x_{q\alpha}\delta_{\gamma\beta}+
x_{p\gamma}x_{q\alpha}\delta_{\tau\beta} \\
& +x_{q\tau}x_{p\beta}\delta_{\gamma\alpha}-x_{q\gamma}x_{p\beta}\delta_{\tau\alpha}-
x_{q\tau}x_{p\alpha}\delta_{\gamma\beta}+x_{q\gamma}x_{p\alpha}\delta_{\tau\beta}),
\end{align*}
where $\epsilon=(-1)^{\alpha+\beta+\gamma+\tau}$.

The above formula is valid for any point $\tilde{\bf x}\in \mathcal{O}(n)$ and not just in the critical points of $G_{\mathcal{O}(n)}$. Even if $\tilde{\bf x}$ is a critical point of $G_{\mathcal{O}(n)}$ it is not necessary a critical point for the prolongation function $G$ and consequently, the partial derivatives of $G$ computed at the point $\tilde{\bf x}$ have a contribution to the formula for the Hessian matrix of the constraint function $G_{\mathcal{O}(n)}$.

\section{Characterization of the critical points for 2-power cost functions defined on ${\bf O}(3)$}

We will exemplify the formulas discovered in the previous section for the case of 2-power cost functions on  ${\bf O}(3)$. 
The orthogonal group ${\bf O}(3)$ is given by:
$${\bf O}(3)=\{X\in \mathcal{M}_{3\times 3}(\R)\,|\,XX^T=\mathbb{I}_3=X^TX\}.$$
A general element of ${\bf O}(3)$ is represented by an orthonormal frame $\{{\bf x}_i=(x_{i1},x_{i2},x_{i3})\}_{i=\overline{1,3}}$ in $\R^3$.
We identify an orthogonal matrix with a vector in $\R^{9}$ using the linear map ${\cal I}:\mathcal{M}_{3\times 3}(\R)\rightarrow \R^{9}$,
\begin{equation}\label{ident}
X\stackrel{{\cal I}}{\longmapsto} \tilde{\bf x}=({\bf x}_1,{\bf x}_2,{\bf x}_3).
\end{equation}
Regarded as a subset of $\R^{9}$, the orthogonal group ${\bf O}(3)$ can be seen as the preimage $\mathcal{O}(3)\subset  \R^{9}$ of the regular value $(\frac{1}{2},\frac{1}{2},\frac{1}{2},0,0,0)\in \R^3\times \R^{3}$ of the constraint functions:
\begin{align*}
F_1 (\tilde{\bf x}) & =\frac{1}{2}||{\bf x}_1||^2,\,\,F_2 (\tilde{\bf x}) =\frac{1}{2}||{\bf x}_2||^2,\,\,F_3 (\tilde{\bf x}) =\frac{1}{2}||{\bf x}_3||^2, \\
F_{12} (\tilde{\bf x}) & = <{\bf x}_1,{\bf x}_2>,\,\,F_{13} (\tilde{\bf x}) = <{\bf x}_1,{\bf x}_3>,\,\,F_{23} (\tilde{\bf x}) = <{\bf x}_2,{\bf x}_3>,
\end{align*}
where $<\cdot,\cdot>$ is the Euclidean product in $\R^3$. The gradients of the constraint functions are given by:
\begin{align*}
\text{\bf grad }F_1(\tilde{\bf x}) & =({\bf x}_1,{\bf 0},{\bf 0}),\,\,\text{\bf grad }F_2(\tilde{\bf x}) =({\bf 0},{\bf x}_2,{\bf 0}),\,\,\text{\bf grad }F_3(\tilde{\bf x})  =({\bf 0},{\bf 0},{\bf x}_3), \\
\text{\bf grad }F_{12}(\tilde{\bf x}) & =({\bf x}_2,{\bf x}_1,{\bf 0}),\,\,\text{\bf grad }F_{13}(\tilde{\bf x}) =({\bf x}_3,{\bf 0},{\bf x}_1),\,\,\text{\bf grad }F_{23}(\tilde{\bf x})  =({\bf 0},{\bf x}_3,{\bf x}_2).
\end{align*}
The Hessian matrices of the constraint functions are given by:
\begin{align*}
\left[\text{Hess}\,{F_1}(\tilde{\bf x})\right] & =
\left[%
\begin{array}{cccc}
  \mathbb{I}_3 & \mathbb{O}_3 & \mathbb{O}_3 \\
  \mathbb{O}_3 & \mathbb{O}_3 & \mathbb{O}_3 \\
 \mathbb{O}_3 & \mathbb{O}_3 & \mathbb{O}_3
\end{array}%
\right],
\left[\text{Hess}\,{F_2}(\tilde{\bf x})\right] =
\left[%
\begin{array}{cccc}
  \mathbb{O}_3 & \mathbb{O}_3 & \mathbb{O}_3 \\
  \mathbb{O}_3 & \mathbb{I}_3 & \mathbb{O}_3 \\
 \mathbb{O}_3 & \mathbb{O}_3 & \mathbb{O}_3
\end{array}%
\right],
\left[\text{Hess}\,{F_3}(\tilde{\bf x})\right] =
\left[%
\begin{array}{cccc}
  \mathbb{O}_3 & \mathbb{O}_3 & \mathbb{O}_3 \\
  \mathbb{O}_3 & \mathbb{O}_3 & \mathbb{O}_3 \\
 \mathbb{O}_3 & \mathbb{O}_3 & \mathbb{I}_3
\end{array}%
\right], \\
\left[\text{Hess}\,{F_{12}}(\tilde{\bf x})\right] & =
\left[%
\begin{array}{cccc}
  \mathbb{O}_3 & \mathbb{I}_3 & \mathbb{O}_3 \\
  \mathbb{I}_3 & \mathbb{O}_3 & \mathbb{O}_3 \\
 \mathbb{O}_3 & \mathbb{O}_3 & \mathbb{O}_3
\end{array}%
\right],
\left[\text{Hess}\,{F_{13}}(\tilde{\bf x})\right] =
\left[%
\begin{array}{cccc}
  \mathbb{O}_3 & \mathbb{O}_3 & \mathbb{I}_3 \\
  \mathbb{O}_3 & \mathbb{O}_3 & \mathbb{O}_3 \\
 \mathbb{I}_3 & \mathbb{O}_3 & \mathbb{O}_3
\end{array}%
\right],
\left[\text{Hess}\,{F_{23}}(\tilde{\bf x})\right] =
\left[%
\begin{array}{cccc}
  \mathbb{O}_3 & \mathbb{O}_3 & \mathbb{O}_3 \\
  \mathbb{O}_3 & \mathbb{O}_3 & \mathbb{I}_3 \\
 \mathbb{O}_3 & \mathbb{I}_3 & \mathbb{O}_3
\end{array}%
\right].
\end{align*}

Considering a cost function  $G_{{\bf O}(3)}:{\bf O}(3)\rightarrow \R$, we identify it with $G_{\mathcal{O}(3)}:=G_{{\bf O}(3)}\circ {\cal I}:\mathcal{O}(3)\rightarrow \R$ and construct a prolongation function $G:\R^9\rightarrow \R$. 
Formula \eqref{Hessian-O(n)} becomes:
\begin{align}
\left[ \text{Hess}\,{G_{\mathcal{O}(3)}}(\tilde{\bf x})\right]  & =\left[ \text{Hess}\,{G}(\tilde{\bf x})\right]_{|_{T_{\tilde{\bf x}}{\mathcal{O}}(3)\times T_{\tilde{\bf x}}{\mathcal{O}}(3)}}-
\sum_{s=1}^3<\frac{\partial G}{\partial {\bf x}_s}(\tilde{\bf x}),{\bf x}_s> \left[ \text{Hess}\,{F_s}(\tilde{\bf x})\right]_{|_{T_{\tilde{\bf x}}{\mathcal{O}}(3)\times T_{\tilde{\bf x}}{\mathcal{O}}(3)}} \nonumber \\
& -\frac{1}{2}\sum_{1\leq p<q\leq 3}(<\frac{\partial G}{\partial {\bf x}_p}(\tilde{\bf x}),{\bf x}_q>+<\frac{\partial G}{\partial {\bf x}_q}(\tilde{\bf x}),{\bf x}_p>) \left[ \text{Hess}\,{F_{pq}}(\tilde{\bf x})\right]_{|_{T_{\tilde{\bf x}}{\mathcal{O}}(3)\times T_{\tilde{\bf x}}{\mathcal{O}}(3)}}.
\end{align}

We choose the base for $T_{\tilde{\bf x}}\mathcal{O}(3)$ as in the previous section:
$$\tilde{\omega}_{12}(\tilde{\bf x})=(x_{12},-x_{11},0,x_{22},-x_{21},0,x_{32},-x_{31},0),$$
 $$\tilde{\omega}_{13}(\tilde{\bf x})=(-x_{13},0,x_{11},-x_{23},0,x_{21},-x_{33},0,x_{31}),$$
$$\tilde{\omega}_{23}(\tilde{\bf x})=(0,x_{13},-x_{12},0,x_{23},-x_{22},0,x_{33},-x_{32}).$$
The restricted Hessian matrices for the constraint functions are the following{\footnote{In the formulas for the restricted Hessian matrices of the constraint functions we have used the fact that $\tilde{\bf x}\in \mathcal{O}(3)$.}}:
$$\left[ \text{Hess}\,{F_1}(\tilde{\bf x})\right]_{|_{T_{\tilde{\bf x}}{\mathcal{O}}(3)\times T_{\tilde{\bf x}}{\mathcal{O}}(3)}}=
\left[%
\begin{array}{cccc}
  1-x_{13}^2 & -x_{12}x_{13} & -x_{11}x_{13} \\
  -x_{12}x_{13} & 1- x_{12}^2 & -x_{11}x_{12} \\
-x_{11}x_{13} &-x_{11}x_{12} &  1-x_{11}^2
\end{array}%
\right],
$$
$$\left[ \text{Hess}\,{F_2}(\tilde{\bf x})\right]_{|_{T_{\tilde{\bf x}}{\mathcal{O}}(3)\times T_{\tilde{\bf x}}{\mathcal{O}}(3)}}=
\left[%
\begin{array}{cccc}
 1- x_{23}^2 & -x_{22}x_{23} & -x_{21}x_{23} \\
  -x_{22}x_{23} &  1-x_{22}^2 & -x_{21}x_{22} \\
-x_{21}x_{23} &-x_{21}x_{22} &  1-x_{21}^2
\end{array}%
\right],
$$
$$\left[ \text{Hess}\,{F_3}(\tilde{\bf x})\right]_{|_{T_{\tilde{\bf x}}{\mathcal{O}}(3)\times T_{\tilde{\bf x}}{\mathcal{O}}(3)}}=
\left[%
\begin{array}{cccc}
  1-x_{33}^2 & -x_{32}x_{33} & -x_{31}x_{33} \\
  -x_{32}x_{33} &  1-x_{32}^2 & -x_{31}x_{32} \\
-x_{31}x_{33} &-x_{31}x_{32} &  1-x_{31}^2
\end{array}%
\right],
$$
$$\left[ \text{Hess}\,{F_{12}}(\tilde{\bf x})\right]_{|_{T_{\tilde{\bf x}}{\mathcal{O}}(3)\times T_{\tilde{\bf x}}{\mathcal{O}}(3)}}=
\left[ \begin {array}{ccc} -2{x_{23}}{x_{13}}&-{x_{23}}{ 
x_{12}}-{x_{13}}{x_{22}}&-{x_{23}}{x_{11}}-{x_{13}}{x_{21}}
\\ -{x_{23}}{x_{12}}-{x_{13}}{x_{22}}&-2{
x_{22}}{x_{12}}&-{x_{22}}{x_{11}}-{x_{12}}{x_{21}}
\\ -{x_{23}}{x_{11}}-{x_{13}}{x_{21}}&-{ 
x_{22}}{x_{11}}-{x_{12}}{x_{21}}&-2{x_{11}}{x_{21}}
\end {array} \right],
$$
$$
\left[ \text{Hess}\,{F_{13}}(\tilde{\bf x})\right]_{|_{T_{\tilde{\bf x}}{\mathcal{O}}(3)\times T_{\tilde{\bf x}}{\mathcal{O}}(3)}}=
\left[ \begin {array}{ccc} -2{x_{33}}{x_{13}}&-{x_{33}}{ 
x_{12}}-{x_{13}}{x_{32}}&-{x_{33}}{x_{11}}-{x_{13}}{x_{31}}
\\ -{x_{33}}{x_{12}}-{x_{13}}{x_{32}}&-2{
x_{32}}{x_{12}}&-{x_{32}}{x_{11}}-{x_{12}}{x_{31}}
\\ -{x_{33}}{x_{11}}-{x_{13}}{x_{31}}&-{ 
x_{32}}{x_{11}}-{x_{12}}{x_{31}}&-2{x_{31}}{x_{11}}
\end {array} \right],
$$
$$
\left[ \text{Hess}\,{F_{23}}(\tilde{\bf x})\right]_{|_{T_{\tilde{\bf x}}{\mathcal{O}}(3)\times T_{\tilde{\bf x}}{\mathcal{O}}(3)}}=
\left[ \begin {array}{ccc} -2{x_{33}}{x_{23}}&-{x_{33}}{
x_{22}}-{x_{23}}{x_{32}}&-{x_{33}}{x_{21}}-{x_{23}}{x_{31}}
\\ -{x_{33}}{x_{22}}-{x_{23}}{x_{32}}&-2{x_{32}}{x_{22}}&-{x_{32}}{x_{21}}-{x_{22}}{x_{31}}
\\ -{x_{33}}{x_{21}}-{x_{23}}{x_{31}}&-{
x_{32}}{x_{21}}-{x_{22}}{x_{31}}&-2{x_{31}}{x_{21}}
\end {array} \right].
$$

We will characterize the critical points of the following 2-power cost function, 
$$G_{\mathbf{O}(3)}(X)=\frac{1}{2}\sum_{i=1}^k ||X-R_i||_F^2,$$
where $R_1,...,R_k$ are sample rotations  and $||\cdot||_F$ is the Frobenius norm. The critical points of the above cost function have been computed using {\it The Embedding Algorithm} in \cite{birtea-comanescu-popa}. Using the identification map $\mathcal{I}$ we obtain the cost function,
$$G_{\mathcal{O}(3)}(\tilde{\bf x})=\frac{1}{2}\sum_{i=1}^k||\tilde{\bf x}-\tilde{\bf r}_i||^2,$$ where $||\cdot||$ is the Euclidean norm on $\R^9$. For the obvious prolongation $G:\R^9\rightarrow \R$ of $G_{\mathcal{O}(3)}$ we have, 
$$\nabla G(\tilde{\bf x})=k(\tilde {\bf x}-\tilde{\bf r}),\,\,
\left[ \text{Hess}\,{G}(\tilde{\bf x})\right]_{|_{T_{\tilde{\bf x}}{\mathcal{O}}(3)\times T_{\tilde{\bf x}}{\mathcal{O}}(3)}}=
\left[ \begin {array}{ccc}
2k & 0 &0  \\
0 & 2k & 0 \\
0 & 0 & 2k
\end {array} \right],
$$
where $\tilde {\bf r}=\frac{1}{k}\sum_{i=1}^k\tilde{\bf r}_i$. 
Applying formula \eqref{Hessian-O(n)} for the case $n=3$, we obtain the components of the Hessian matrix of the cost function $G_{\mathcal{O}(3)}$:
\begin{align*}
h_{11}(\tilde{\bf x}) & =k(x_{11}r_{11}+x_{21}r_{21}+x_{31}r_{31}+x_{12}r_{12}+x_{22}r_{22}+x_{32}r_{32}), \\
h_{12}(\tilde{\bf x}) & =-\frac{k}{2}(x_{12}r_{13}+x_{22}r_{23}+x_{32}r_{33}+x_{13}r_{12}+x_{23}r_{22}+x_{33}r_{32}), \\
h_{13}(\tilde{\bf x}) & =-\frac{k}{2}(x_{11}r_{13}+x_{21}r_{23}+x_{31}r_{33}+x_{13}r_{11}+x_{23}r_{21}+x_{33}r_{31}), \\
h_{22}(\tilde{\bf x}) & =k(x_{11}r_{11}+x_{21}r_{21}+x_{31}r_{31}+x_{13}r_{13}+x_{23}r_{23}+x_{33}r_{33}), \\
h_{23}(\tilde{\bf x}) & =-\frac{k}{2}(x_{12}r_{11}+x_{22}r_{21}+x_{32}r_{31}+x_{11}r_{12}+x_{21}r_{22}+x_{31}r_{32}), \\
h_{33}(\tilde{\bf x}) & =k(x_{12}r_{12}+x_{22}r_{22}+x_{32}r_{32}+x_{13}r_{13}+x_{23}r_{23}+x_{33}r_{33}).
\end{align*}

For the columns of the matrices $X$, respectively $R=\frac{1}{k}\sum_{i=1}^kR_i$, we make the notations ${\bf y}_i=(x_{1i},x_{2i},x_{3i})$, and respectively  ${\bf s}_i=(r_{1i},r_{2i},r_{3i})$. The components of the of the Hessian matrix of the cost function $G_{\mathcal{O}(3)}$ can be written in the equivalent form:
\begin{align*}
h_{11}(\tilde{\bf x}) & =k(<{\bf y}_1,{\bf s}_1>+<{\bf y}_2,{\bf s}_2>), \\
h_{12}(\tilde{\bf x}) & =-\frac{k}{2}(<{\bf y}_2,{\bf s}_3>+<{\bf y}_3,{\bf s}_2>), \\
h_{13}(\tilde{\bf x}) & =-\frac{k}{2}(<{\bf y}_1,{\bf s}_3>+<{\bf y}_3,{\bf s}_1>), \\
h_{22}(\tilde{\bf x}) & =k(<{\bf y}_1,{\bf s}_1>+<{\bf y}_3,{\bf s}_3>), \\
h_{23}(\tilde{\bf x}) & =-\frac{k}{2}(<{\bf y}_2,{\bf s}_1>+<{\bf y}_1,{\bf s}_2>), \\
h_{33}(\tilde{\bf x}) & =k(<{\bf y}_2,{\bf s}_2>+<{\bf y}_3,{\bf s}_3>).
\end{align*}

\begin{rem}
The above expressions for the Hessian matrix depend on the chosen base for the tangent space $T_{\tilde{\bf x}}\mathcal{O}(3)$. If we rename the base chosen above as follows:
$$\boldsymbol{\nu}_1(\tilde{\bf x})=\tilde{\omega}_{23}(\tilde{\bf x}),\,\,\boldsymbol{\nu}_2(\tilde{\bf x})=\tilde{\omega}_{13}(\tilde{\bf x}),\,\,\boldsymbol{\nu}_3(\tilde{\bf x})=\tilde{\omega}_{12}(\tilde{\bf x}),$$
the formulas for the  components of the Hessian matrix of the cost function $G_{\mathcal{O}(3)}$ have a more natural expressions with respect to the symmetry of $\mathbf{O}(3)$: 
\begin{align*}
\overline{h}_{11}(\tilde{\bf x}) & =k(<{\bf y}_2,{\bf s}_2>+<{\bf y}_3,{\bf s}_3>), \nonumber \\
\overline{h}_{12}(\tilde{\bf x}) & =-\frac{k}{2}(<{\bf y}_1,{\bf s}_2>+<{\bf y}_2,{\bf s}_1>), \nonumber \\
\overline{h}_{13}(\tilde{\bf x}) & =-\frac{k}{2}(<{\bf y}_1,{\bf s}_3>+<{\bf y}_3,{\bf s}_1>), \nonumber \\
\overline{h}_{22}(\tilde{\bf x}) & =k(<{\bf y}_1,{\bf s}_1>+<{\bf y}_3,{\bf s}_3>), \nonumber \\
\overline{h}_{23}(\tilde{\bf x}) & =-\frac{k}{2}(<{\bf y}_2,{\bf s}_3>+<{\bf y}_3,{\bf s}_2>), \nonumber \\
\overline{h}_{33}(\tilde{\bf x}) & =k(<{\bf y}_1,{\bf s}_1>+<{\bf y}_2,{\bf s}_2>). \Box
\end{align*} 
\end{rem}
Using the intrinsic Riemannian geometry of the Lie group ${\bf SO}(3)$, another formula for the Hessian matrix of the 2-power cost function has been given in \cite{moakher}.
\medskip

Now we apply the above formulas to determine the nature of the critical points of the following example of 2-power cost function defined on the connected component of the identity matrix of the orthogonal group ${\bf O}(3)$ which is ${\bf SO}(3)$:
$$G_{\mathbf{SO}(3)}^{\alpha}(X)=\frac{1}{2}(||X-R_1||_F^2+||X-R_2||_F^2+||X-R_3||_F^2),$$
where 
\begin{equation*}
R_1=\left(
\begin{array}{ccc}
1 & 0 & 0 \\
0  & -1 & 0 \\
0 & 0 & -1
\end{array}
\right),\,\,
R_2=\left(
\begin{array}{ccc}
1 & 0 & 0 \\
0  & 0 & -1 \\
0 & 1 & 0
\end{array}
\right),\,\,
R_3=\left(
\begin{array}{ccc}
1 & 0 & 0 \\
0  & \cos\alpha & -\sin\alpha \\
0 & \sin\alpha & \cos \alpha
\end{array}
\right),\,\,\alpha\in [-\pi,\pi],
\end{equation*}
are rotations along the $x$-axis.  Using the identification map $\mathcal{I}$, we obtain the cost function defined on the connected component $\mathcal{SO}(3)$ of the point $(1,0,0,0,1,0,0,0,1)\in\mathcal{O}(3)$:
$$G_{\mathcal{SO}(3)}^{\alpha}(\tilde{\bf x})=\frac{1}{2}(||\tilde{\bf x}-\tilde{\bf r}_1||^2+||\tilde{\bf x}-\tilde{\bf r}_2||^2+||\tilde{\bf x}-\tilde{\bf r}_3||^2),$$ where $||\cdot||$ is the Euclidean norm on $\R^9$. For the obvious prolongation $G^{\alpha}:\R^9\rightarrow \R$ of $G_{\mathcal{SO}(3)}^{\alpha}$ we have:
$$\nabla G^{\alpha}(\tilde{\bf x})=3(\tilde {\bf x}-\tilde{\bf r}),\,\,
\left[ \text{Hess}\,{G^{\alpha}}(\tilde{\bf x})\right]_{|_{T_{\tilde{\bf x}}{\mathcal{O}}(3)\times T_{\tilde{\bf x}}{\mathcal{O}}(3)}}=
\left[ \begin {array}{ccc}
6 & 0 &0  \\
0 & 6 & 0 \\
0 & 0 & 6
\end {array} \right],
$$
where $\tilde {\bf r}=\frac{1}{3}(\tilde{\bf r}_1+\tilde{\bf r}_2+\tilde{\bf r}_3)=(1,0,0,0,\frac{-1+\cos\alpha}{3},\frac{-1-\sin\alpha}{3},0,\frac{1+\sin\alpha}{3},\frac{-1+\cos\alpha}{3})$. 

By using the base $\{\boldsymbol{\nu}_1(\tilde{\bf x}),\boldsymbol{\nu}_2(\tilde{\bf x}),\boldsymbol{\nu}_3(\tilde{\bf x})\}$ for the tangent space $T_{\tilde{\bf x}}\mathcal{SO}(3)$, we obtain the coefficients of the Hessian matrix of the cost function $G_{\mathcal{SO}(3)}^{\alpha}$:
\begin{align}\label{Hess-baza-buna}
\overline{h}_{11}(\tilde{\bf x}) & =(-1+\cos\alpha)(x_{22}+x_{33})+(1+\sin\alpha)(x_{32}-x_{23}), \nonumber \\
\overline{h}_{12}(\tilde{\bf x}) & =-\frac{1}{2}\left[(1+\sin\alpha)x_{31}+(-1+\cos\alpha)x_{21}+3x_{12}\right], \nonumber \\
\overline{h}_{13}(\tilde{\bf x}) & =-\frac{1}{2}\left[(-1+\cos\alpha)x_{31}-(1+\sin\alpha)x_{21}+3x_{13}\right],  \\
\overline{h}_{22}(\tilde{\bf x}) & =(-1+\cos\alpha)x_{33}-(1+\sin\alpha)x_{23}+3x_{11}, \nonumber \\
\overline{h}_{23}(\tilde{\bf x}) & =-\frac{1}{2}\left[(-1+\cos\alpha)(x_{32}+x_{23})+(1+\sin\alpha)(x_{33}-x_{22})\right], \nonumber \\
\overline{h}_{33}(\tilde{\bf x}) & =(1+\sin\alpha)x_{32}+(-1+\cos\alpha)x_{22}+3x_{11}. \nonumber
\end{align}

The critical points of the 2-power cost function $G_{\mathbf{SO}(3)}^{\alpha}$ have been computed in \cite{birtea-comanescu-popa} using the Embedding Algorithm.  We find five sets of critical points as follows:
\begin{align*}
& \text{Set }_{black}^{Rot} = \{R^{\bf q}\,|\,{\bf q}=(0,0,\pm\sqrt{1-t^2},t),\,t\in[-1,1]\}, \\
& \text{Set }_{green}^{Rot} = \{R^{\bf q}\,|\,{\bf q}=(\sqrt{1-x_{2,min}^2(\alpha)},x_{2,min}(\alpha),0,0),\,\alpha\in[-\pi,\pi]\}, \\
& \text{Set }_{pink}^{Rot} = \{R^{\bf q}\,|\,{\bf q}=(-\sqrt{1-x_{2,min}^2(\alpha)},x_{2,min}(\alpha),0,0),\,\alpha\in[-\pi,\pi]\}, \\
& \text{Set}_{red}^{Rot} = \{R^{\bf q}\,|\,{\bf q}=(\sqrt{1-x_{2,max}^2(\alpha)},x_{2,max}(\alpha),0,0),\,\alpha\in[-\pi,\pi]\}, \\
& \text{Set}_{blue}^{Rot} = \{R^{\bf q}\,|\,{\bf q}=(-\sqrt{1-x_{2,max}^2(\alpha)},x_{2,max}(\alpha),0,0),\,\alpha\in[-\pi,\pi]\},
\end{align*}
where $x_{2,min}(\alpha)$ and $x_{2,max}(\alpha)$ are the smallest, respectively largest real positive solutions of the polynomial
\begin{align*}
Q_{2,\alpha}(Z) & =\left( 128 \sin^4 \frac{\alpha}{2}-32
 \sin^2 \frac{\alpha}{2}+4 \right) Z^{4}-\left( 128 \sin^4 \frac{\alpha}{2}-32
 \sin^2 \frac{\alpha}{2}+4 \right) Z^{2} \\
& -16\sin^6 \frac{\alpha}{2}+16\sin^5 \frac{\alpha}{2}\cos \frac{\alpha}{2}+28\sin^4 \frac{\alpha}{2}-8 \sin^2 \frac{\alpha}{2} +1,
\end{align*}
and $R^{\bf q}$ is the rotation corresponding to the unit quaternion ${\bf q}$, see formula \eqref{Rq} in Annexe. 

We study the nature of the above critical points using the Hessian characterization.
\medskip

{\bf Case black.} The critical points corresponding to rotations from $ \text{Set }_{black}^{Rot}$ are absolute maximum for the cost function $G_{\mathbf{SO}(3)}^{\alpha}$ as have been pointed out in \cite{birtea-comanescu-popa}. Applying formula \eqref{Hess-baza-buna} at critical points in $ \mathcal{I}(\text{Set }_{black}^{Rot})$ we obtain the eigenvalues $\lambda_1=0$, $\lambda_2=-3+\sqrt{3+2\sin\alpha-2\cos\alpha}$, $\lambda_3=-3-\sqrt{3+2\sin\alpha-2\cos\alpha}$. Consequently, a critical point of this set is a degenerate absolute maximum. 
\begin{figure}[H]
\begin{center}
  \includegraphics[width=15cm]{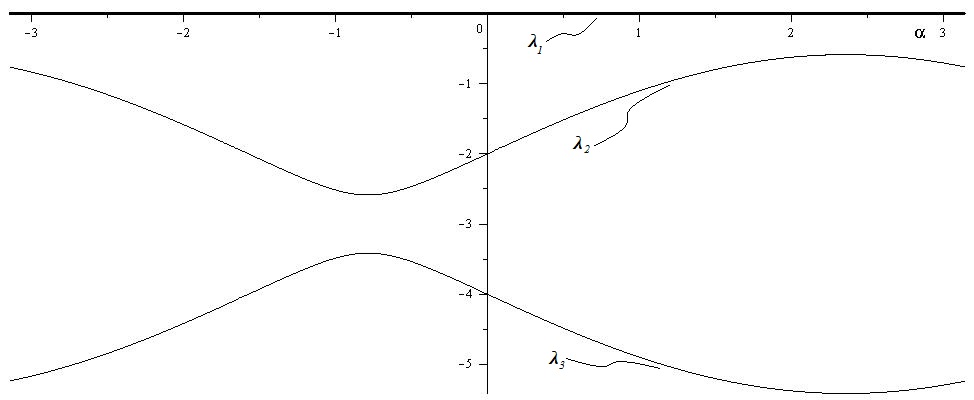}\\
  \caption{Eigenvalues of the Hessian matrix of $G_{\mathbf{SO}(3)}^{\alpha}$ computed at critical points from $ \text{Set }_{black}^{Rot}$.}\label{Set-black}
\end{center}
\end{figure}
\medskip

\begin{figure}[H]
\begin{center}
  \includegraphics[width=15cm]{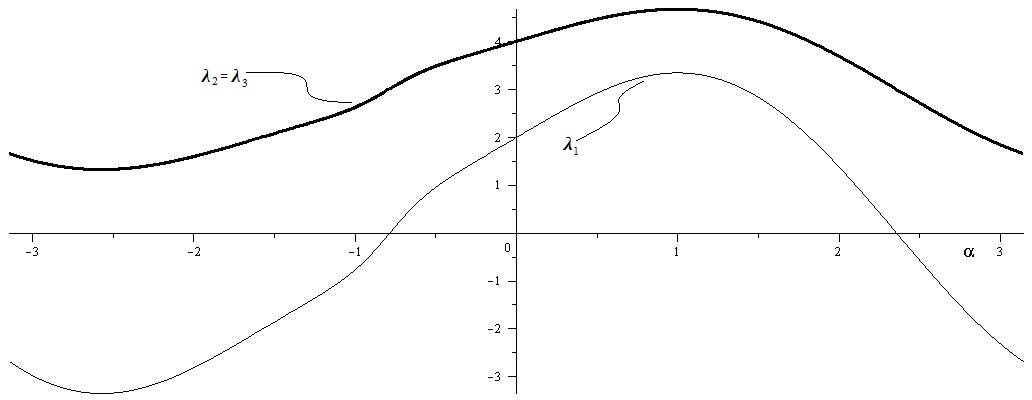}\\
  \caption{Eigenvalues of the Hessian matrix of $G_{\mathbf{SO}(3)}^{\alpha}$ computed at critical points from $ \text{Set }_{green}^{Rot}$.}\label{Set-green}
\end{center}
\end{figure}

{\bf Case green.} 
Applying formula \eqref{Hess-baza-buna} at critical points in $ \mathcal{I}(\text{Set }_{green}^{Rot})$ we obtain for the Hessian matrix of the cost function  $G_{\mathbf{SO}(3)}^{\alpha}$ two equal eigenvalues that are represented by the thick line in the Figure \ref{Set-green} and one simple eigenvalue. 
We note that for $\alpha=-\frac{\pi}{4}$ and $\alpha=\frac{3\pi}{4}$ we obtain a bifurcation phenomena and the critical points in  $ \text{Set }_{green}^{Rot}$ corresponding to this values of the parameter $\alpha$ have a degenerate Hessian matrix. For $\alpha \in [-\pi,-\frac{\pi}{4})\cup (\frac{3\pi}{4},\pi]$ the corresponding critical points in  $ \text{Set }_{green}^{Rot}$ are saddle critical points for the cost function $G_{\mathbf{SO}(3)}^{\alpha}$.  For $\alpha \in (-\frac{\pi}{4},\frac{3\pi}{4})$ the corresponding critical points are local minima.

{\bf Case pink.} For the critical points in the set $ \text{Set }_{pink}^{Rot}$ the Hessian matrix of the cost function $G_{\mathbf{SO}(3)}^{\alpha}$ has one negative eigenvalue and two equal positive eigenvalues. Consequently, this critical points are all saddle points.
\begin{figure}[H]
\begin{center}
  \includegraphics[width=15cm]{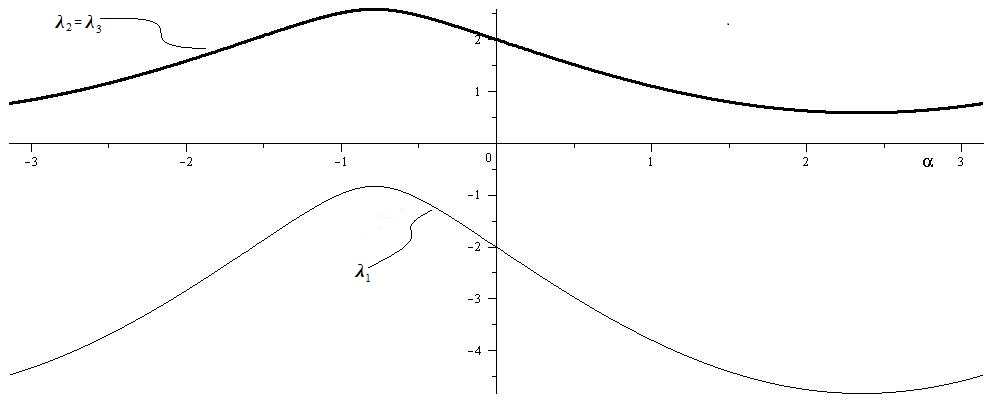}\\
  \caption{Eigenvalues of the Hessian matrix of $G_{\mathbf{SO}(3)}^{\alpha}$ computed at critical points from $ \text{Set }_{pink}^{Rot}$.}\label{Set-pink}
\end{center}
\end{figure}

{\bf Case red.} The critical points in the  set $ \text{Set }_{red}^{Rot}$ are all local minima as eigenvalues of the Hessian matrix of the cost function $G_{\mathbf{SO}(3)}^{\alpha}$ are all positive. 
\begin{figure}[H]
\begin{center}
  \includegraphics[width=15cm]{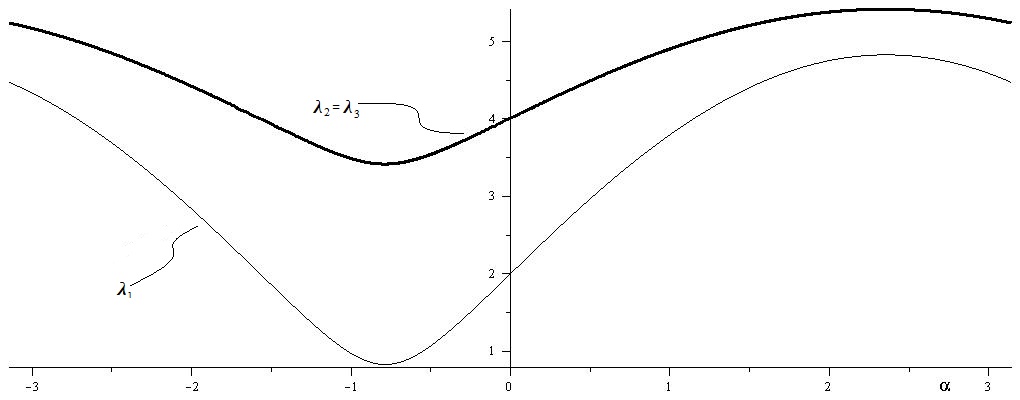}\\
  \caption{Eigenvalues of the Hessian matrix of $G_{\mathbf{SO}(3)}^{\alpha}$ computed at critical points from $ \text{Set }_{red}^{Rot}$.}\label{Set-red}
\end{center}
\end{figure}

{\bf Case blue.} Again a bifurcation phenomena appears for this case at the values $\alpha=-\frac{\pi}{4}$ and $\alpha=\frac{3\pi}{4}$.  For $\alpha \in [-\pi,-\frac{\pi}{4})\cup (\frac{3\pi}{4},\pi]$ the corresponding critical points in  $ \text{Set }_{blue}^{Rot}$ are local minima for the cost function $G_{\mathbf{SO}(3)}^{\alpha}$.  For $\alpha \in (-\frac{\pi}{4},\frac{3\pi}{4})$ the corresponding critical points are saddle points.
\begin{figure}[H]
\begin{center}
  \includegraphics[width=15cm]{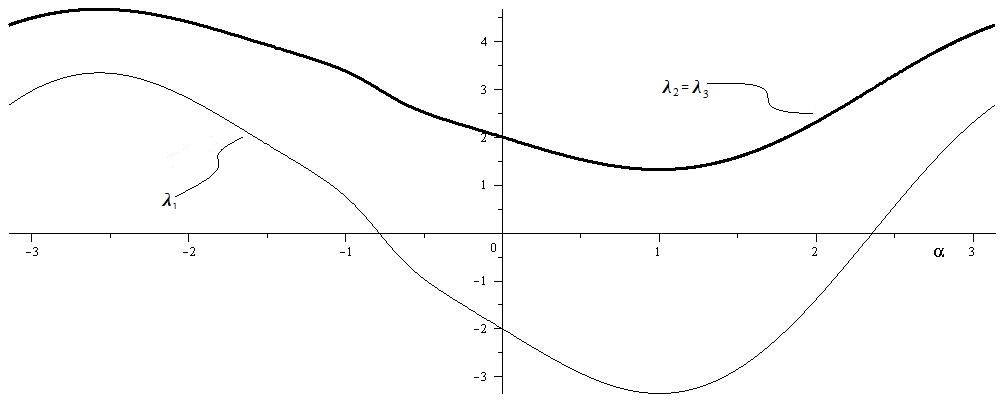}\\
  \caption{Eigenvalues of the Hessian matrix of $G_{\mathbf{SO}(3)}^{\alpha}$ computed at critical points from $ \text{Set }_{blue}^{Rot}$.}\label{Set-blue}
\end{center}
\end{figure}

\section{Stability of equilibrium points using restricted Hessian}

We apply the results of Section 2 to the stability problem of an equilibrium point for a dynamical system generated by a vector field $X_S$ defined on a manifold $S$.
Let $x_e\in S$ be an equilibrium point for the dynamics on the manifold $S$ generated by the vector field $X_S$. Stability behavior of the equilibrium point $x_e$ can be determined using the direct method of Lyapunov. This method requires the knowledge of a Lyapunov function $G_S:S\rightarrow \R$ which has the following properties:
\begin{itemize}
\item [(i)] $\dot{G}_S:={\bf L}_{X_S}G_S\leq 0$;

\item [(ii)] $G_S(x)>G_S(x_e)$, for all $x\neq x_e$ in a neighborhood of $x_e$.
\end{itemize}
In order to verify the above conditions one needs to construct a local system of coordinates on $S$ around the equilibrium point $x_e$. If $G_S$ is a $C^2$ differentiable function and $x_e$ is a critical point for $G_S$ then a sufficient condition for $(\text{ii})$ to hold is given by the positive definiteness of the Hessian matrix $\left[\text{Hess}\,{G_S}(x_e)\right]$. Usually is very difficult to construct local coordinates on the submanifold $S$ and in these cases we will bypass this difficulty by embedding the problem in an ambient space $M$ (usually an Euclidean space) and use formula  \eqref{Hessian-constraints} given in Theorem \ref{calcul-hessiana-restransa}. 

 Suppose that the manifold $S$ is a preimage of a regular value for a map ${\bf F}=(F_1,...,F_k):M\rightarrow \R^k$, where $(M,g)$ is an ambient Riemannian manifold. Let $X\in \mathcal{X}(M)$ be a prolongation of the vector field $X_S$, i.e. $X_{|S}=X_S\in \mathcal{X}(S)$ and $G:M\rightarrow \R$ be a $C^2$ prolongation of the function $G_S$. The equilibrium point $x_e$ is also an equilibrium point for the dynamics on $M$ generated by the vector field $X$. 
The conditions of the next result guaranties the applicability of the direct method of Lyapunov stated above.

\begin{thm}\label{stability}
The following conditions:
\begin{itemize}
\item [(i)] $\dot{G}:={\bf L}_XG\leq 0$,

\item [(ii)] $\text{\bf grad }G(x_e)=\sum_{i=1}^k\sigma_i(x_e)\text{\bf grad }F_i(x_e)$,

\item [(iii)] $\left[ \normalfont{\text{Hess}}\,{G}(x_e)\right]_{|_{T_{_{x_e}}S\times T_{_{x_e}}S}}-\sum_{i=1}^k\sigma_i(x_e)\left[ \normalfont{\text{Hess}}\,{F_i}(x_e)\right]_{|_{T_{_{x_e}}S\times T_{_{x_e}}S}}$ is positive definite,
\end{itemize}
implies that the equilibrium point $x_e$  is stable for the dynamics generated by the vector field $X_S$. 
\end{thm}

The condition ${\bf L}_{X_S}G_S\leq 0$ is implied by the condition $(i)$ in the above theorem. Condition $(ii)$ is equivalent with $x_e$ being a critical point of the function $G_S:S\rightarrow \R$, where $\sigma_i(x_e)$ are the Lagrange multipliers, and condition $(iii)$ is equivalent with positive definiteness of the Hessian matrix $\left[\text{Hess}\,{G_S}(x_e)\right]$.
The advantage of the above theorem is that all the necessary computations for verifying conditions $(i)$, $(ii)$ and $(iii)$ are made using the coordinates of the ambient space $M$ which usually is an Euclidean space. {\bf Note that the constraint functions $F_1,...,F_k$ do not need to be conserved quantities for the prolonged vector field $X$.}
\medskip

Usually the above theorem is applied backwards, where the vector field $X_S$ is the restriction of a vector field $X\in \mathcal{X}(M)$ to an invariant submanifold $S$ under the dynamics generated by the vector field $X$.
In the case when $F_1,...,F_k,G$ are conserved quantities for the vector field $X$ and the conditions $(ii)$ and $(iii)$ of the above theorem are satisfied, then the equilibrium point $x_e$ is also stable for the dynamics generated by the vector field $X$ according to the algebraic method, see \cite{comanescu}, \cite{comanescu-1}, \cite{comanescu-2}. 
\medskip

We will apply the above result to the following Hamilton-Poisson situation. Let $(M,\{\cdot,\cdot\})$ be a finite dimensional Poisson manifold and $X_H$ a Hamilton-Poisson vector field. The paracompactness of the manifold $M$ ensures the existence of a Riemannian metric. The conserved quantities are the Casimir functions $F_1=C_1,...,F_k=C_k$. A regular symplectic leaf $S$ is an open dense set of the submanifold generated by a regular value of the Casimir functions (usually the two sets are equal). The restricted vector field $X_S=(X_H)_{_{|S}}$ on the symplectic leaf is again a Hamiltonian vector field with respect to the symplectic form induced by the Poisson structure $\{\cdot,\cdot\}$ and the Hamiltonian function $G_S=H_{|S}$.
If $x_e\in S$ is an equilibrium point for the restricted Hamiltonian vector field $X_S$, then it is also an equilibrium point for $X_H$ and also it is a critical point of the restricted Hamiltonian function $G_S=H_{|S}$. 
We are in the hypotheses of the Theorem \ref{stability} and the algebraic method for stability.

\begin{thm}
A sufficient condition for stability of the equilibrium point $x_e$ with respect to the dynamics $X_H$ is given by the following condition:
$$\left(\left[ \normalfont{\text{Hess}}\,{H}(x_e)\right]-\sum_{i=1}^k\sigma_i(x_e)\left[ \normalfont{\text{Hess}}\,{C_i}(x_e)\right]\right)_{|_{T_{_{x_e}}S\times T_{_{x_e}}S}}\,\, \text{is positive definite}.$$
\end{thm}

The matrix in the above theorem represents the Hessian matrix of the restricted function $H_{|S}:S\rightarrow \R$ as has been shown in Theorem \ref{calcul-hessiana-restransa}. It is also the Hessian matrix restricted to the tangent space $ T_{_{x_e}}S$ of the augmented function $F:M\rightarrow \R$, $F(x)=H(x)-\sum_{i=1}^k\sigma_i(x_e)C_i(x)$ used in \cite{maddocks-1991}, \cite{beck-hall}, and \cite{wang-xu}. 

In the case of symmetries the augmented function method for stability of relative equilibria has been studied extensively in \cite{patrick}, \cite{ortega-1}, \cite{ortega-2}, \cite{patrick-roberts}, \cite{olmos}, \cite{montaldi}. The tangent space to the invariant submanifold $S$ can be further decomposed taking into account the symmetry of the dynamical system under study.

\section{Annexe}
{\bf I.} The Gramian associated to the constraint functions \eqref{constraints-1}-\eqref{constraints-2} is  $\Sigma(\tilde{\bf x}):=\det\Sigma_{(F_1,\ldots ,F_n,F_{12},...,F_{n-1,n})}^{(F_1,\ldots ,F_n,F_{12},...,F_{n-1,n})}(\tilde{\bf x})$. The matrix  $\Sigma_{(F_1,\ldots ,F_n,F_{12},...,F_{n-1,n})}^{(F_1,\ldots ,F_n,F_{12},...,F_{n-1,n})}(\tilde{\bf x})$ has the form
$$\Sigma_{(F_1,\ldots ,F_n,F_{12},...,F_{n-1,n})}^{(F_1,\ldots ,F_n,F_{12},...,F_{n-1,n})}(\tilde{\bf x})=
\left[%
\begin{array}{cc}
 A & C^T \\
 C & B  \\
\end{array}%
\right],$$
where
\begin{align*}
A & =\left[
\begin{array}{ccc}
 <\text{\bf grad }F_{1},\text{\bf grad }F_{1}> & ... & <\text{\bf grad }F_{n},\text{\bf grad }F_{1}> \\
...& ... & ... \\
 <\text{\bf grad }F_{1},\text{\bf grad }F_{n}> & ... & <\text{\bf grad }F_{n},\text{\bf grad }F_{n}> \\
\end{array}\right], \\
B & =\left[
\begin{array}{ccc}
 <\text{\bf grad }F_{12},\text{\bf grad }F_{12}> & ... & <\text{\bf grad }F_{n-1,n},\text{\bf grad }F_{12}> \\
...& ... & ... \\
 <\text{\bf grad }F_{12},\text{\bf grad }F_{n-1,n}> & ... & <\text{\bf grad }F_{n-1,n},\text{\bf grad }F_{n-1,n}> \\
\end{array}\right], \\
C & =\left[
\begin{array}{ccc}
 <\text{\bf grad }F_{1},\text{\bf grad }F_{12}> & ... & <\text{\bf grad }F_{n},\text{\bf grad }F_{12}> \\
...& ... & ... \\
 <\text{\bf grad }F_{1},\text{\bf grad }F_{n-1,n}> & ... & <\text{\bf grad }F_{n},\text{\bf grad }F_{n-1,n}> \\
\end{array}\right]. \\
\end{align*}
We have the following computations:
\begin{align*}
 <\text{\bf grad }F_{s}(\tilde{\bf x}),\text{\bf grad }F_{r}(\tilde{\bf x})> & =<\sum_{i=1}^nx_{si}  \tilde{\bf e}_{si},\sum_{j=1}^nx_{rj}  \tilde{\bf e}_{rj}> 
=\sum_{i,j=1}^nx_{si}x_{rj}< \tilde{\bf e}_{si}, \tilde{\bf e}_{rj}> \\
& =\sum_{i,j=1}^nx_{si}x_{rj}\delta_{sr}\delta_{ij}=\sum_{i=1}^nx_{si}x_{ri}\delta_{sr} \\
& = <{\bf x}_s,{\bf x}_r>\delta_{sr},
\end{align*}
\begin{align*}
 <\text{\bf grad }F_{s}(\tilde{\bf x}),\text{\bf grad }F_{\alpha\beta}(\tilde{\bf x})> & =<\sum_{i=1}^nx_{si}  \tilde{\bf e}_{si},\sum_{j=1}^n(x_{\beta j}  \tilde{\bf e}_{\alpha j}+x_{\alpha j}  \tilde{\bf e}_{\beta j})> \\
& =\sum_{i=1}^n(x_{si}x_{\beta i}\delta_{s\alpha}+x_{si}x_{\alpha i}\delta_{s\beta}) \\
& = <{\bf x}_s,{\bf x}_{\beta}>\delta_{s\alpha}+ <{\bf x}_s,{\bf x}_{\alpha}>\delta_{s\beta},
\end{align*}
\begin{align*}
 <\text{\bf grad }F_{\gamma\tau}(\tilde{\bf x}),\text{\bf grad }F_{\alpha\beta}(\tilde{\bf x})> & =<\sum_{i=1}^n(x_{\tau i}  \tilde{\bf e}_{\gamma i}+x_{\gamma i}  \tilde{\bf e}_{\tau i}),\sum_{j=1}^n(x_{\beta j}  \tilde{\bf e}_{\alpha j}+x_{\alpha j}  \tilde{\bf e}_{\beta j})> \\
& =\sum_{i=1}^n(x_{\tau i}x_{\beta i}\delta_{\gamma\alpha}+x_{\tau i}x_{\alpha i}\delta_{\gamma\beta}+x_{\gamma i}x_{\beta i}\delta_{\tau\alpha}+x_{\gamma i}x_{\alpha i}\delta_{\tau\beta}) \\
& = <{\bf x}_{\tau},{\bf x}_{\beta}>\delta_{\gamma\alpha}+ <{\bf x}_{\tau},{\bf x}_{\alpha}>\delta_{\gamma\beta}+<{\bf x}_{\gamma},{\bf x}_{\beta}>\delta_{\tau\alpha}+ <{\bf x}_{\gamma},{\bf x}_{\alpha}>\delta_{\tau\beta}.
\end{align*}

Consequently, using the identification \eqref{ident} for $\tilde{\bf x}\in \mathcal{O}(n)={\cal I}({\bf O}(n))$ we have,
$$\Sigma(\tilde{\bf x})=\det\left[%
\begin{array}{cc}
 \mathbb{I}_n & \mathbb{O} \\
 \mathbb{O} & 2\mathbb{I}_{\frac{n(n-1)}{2}}  \\
\end{array}%
\right]=2^{\frac{n(n-1)}{2}}.
$$

For $\tilde{\bf x}\in \mathcal{O}(n)$ we have the following computations:
\begin{align*}
\Sigma_s(\tilde{\bf x}) & =\det\Sigma_{(F_1,...,F_{s-1},G,F_{s+1},... ,F_n,F_{12},...,F_{n-1,n})}^{(F_1,\ldots ,F_n,F_{12},...,F_{n-1,n})}(\tilde{\bf x}) \\
& =  
\det
\left[
\begin{array}{c|c}
\begin{array}{ccccc}
1 & ... &<\text{\bf grad }G(\tilde{\bf x}),\text{\bf grad }F_{1}(\tilde{\bf x})> & ...& 0 \\
...& ... & ... & ... & ... \\
 0 & ... &<\text{\bf grad }G(\tilde{\bf x}),\text{\bf grad }F_{s}(\tilde{\bf x})> & ...& 0 \\
...& ... & ... & ... & ...\\
 0 & ... &<\text{\bf grad }G(\tilde{\bf x}),\text{\bf grad }F_{n}(\tilde{\bf x})> & ...& 1 
\end{array} & 
\begin{array}{ccc}
0 & ... & 0 \\
... & ... & ... \\
0 & ... & 0 \\
... & ... & ... \\
0 & ... & 0
\end{array} \\
\hline \\
\begin{array}{ccccc}
0 & ... &<\text{\bf grad }G(\tilde{\bf x}),\text{\bf grad }F_{12}(\tilde{\bf x})> & ...& 0 \\
...& ... & ... & ... & ... \\
 0 & ... &<\text{\bf grad }G(\tilde{\bf x}),\text{\bf grad }F_{n-1,n}(\tilde{\bf x})> & ...& 0 
\end{array} & 
\begin{array}{ccc}
2 & ... & 0 \\
... & ... & ... \\
0 & ... & 2
\end{array} \\
\end{array}\right] \\
& = <\text{\bf grad }G(\tilde{\bf x}),\text{\bf grad }F_{s}(\tilde{\bf x})> \det\left[
\begin{array}{cc}
 \mathbb{I}_{n-1} & \mathbb{O} \\
 \mathbb{O} & 2\mathbb{I}_{\frac{n(n-1)}{2}}  \\
\end{array} \right] \\
& =2^{\frac{n(n-1)}{2}}<\text{\bf grad }G(\tilde{\bf x}),\text{\bf grad }F_{s}(\tilde{\bf x})> ,
\end{align*}
\begin{align*}
\Sigma_{pq}(\tilde{\bf x}) & =\det\Sigma_{(F_1,... ,F_n,F_{12},...,G,...,F_{n-1,n})}^{(F_1,\ldots ,F_n,F_{12},...,F_{n-1,n})}(\tilde{\bf x}) \\
& = <\text{\bf grad }G(\tilde{\bf x}),\text{\bf grad }F_{pq}(\tilde{\bf x})> \det\left[
\begin{array}{cc}
 \mathbb{I}_{n} & \mathbb{O} \\
 \mathbb{O} & 2\mathbb{I}_{\frac{n(n-1)}{2}-1}  \\
\end{array} \right] \\
& =2^{\frac{n(n-1)}{2}-1}<\text{\bf grad }G(\tilde{\bf x}),\text{\bf grad }F_{pq}(\tilde{\bf x})> .
\end{align*}
Consequently, 
\begin{equation}\label{sigma-uri}
\sigma_s(\tilde{\bf x})= <\text{\bf grad }G(\tilde{\bf x}),\text{\bf grad }F_{s}(\tilde{\bf x})> ,\,\,
\sigma_{pq}(\tilde{\bf x})=\frac{1}{2} <\text{\bf grad }G(\tilde{\bf x}),\text{\bf grad }F_{pq}(\tilde{\bf x})>.
\end{equation}

{\bf II.} We have the following formula for the multiplication of the two $n\times n$ matrices ${\bf e}_i\otimes{\bf e}_j$ and ${\bf e}_{\alpha}\otimes{\bf e}_{\beta}$:
\begin{equation*}\label{mutiplicare-matrici-baza}
{\bf e}_i\otimes{\bf e}_j\cdot {\bf e}_{\alpha}\otimes{\bf e}_{\beta}=\delta_{j\alpha}{\bf e}_i\otimes{\bf e}_{\beta}.
\end{equation*}

Using \eqref{Omega-100}, a base for $T_X{\bf O}(n)$ is given by he following matrices,
\begin{align}\label{calcul-baza-tangenta}
X\Omega_{\alpha\beta} & = (\sum_{i,j=1}^nx_{ij}{\bf e}_i\otimes {\bf e}_j)((-1)^{\alpha+\beta}({\bf e}_{\alpha}\otimes {\bf e}_{\beta}-{\bf e}_{\beta}\otimes {\bf e}_{\alpha})) \nonumber \\
& = (-1)^{\alpha+\beta}\sum_{i,j=1}^nx_{ij}({\bf e}_i\otimes{\bf e}_j\cdot {\bf e}_{\alpha}\otimes{\bf e}_{\beta}-{\bf e}_i\otimes{\bf e}_j\cdot {\bf e}_{\beta}\otimes{\bf e}_{\alpha}) \nonumber \\
& = (-1)^{\alpha+\beta}\sum_{i,j=1}^nx_{ij}(\delta_{j\alpha}{\bf e}_i\otimes{\bf e}_{\beta}-\delta_{j\beta}{\bf e}_i\otimes{\bf e}_{\alpha}) \nonumber \\
& = (-1)^{\alpha+\beta}\sum_{i=1}^n(x_{i\alpha}{\bf e}_i\otimes{\bf e}_{\beta}-x_{i\beta}{\bf e}_i\otimes{\bf e}_{\alpha}).
\end{align}

In order to compute the restricted Hessian to the tangent space $T_{\tilde{\bf x}}{\mathcal O}(n)$ we need the following computation,
\begin{align*}
<\tilde{\omega}_{\gamma\tau}(\tilde{\bf x}), & \tilde{\bf e}_{ab}\otimes \tilde{\bf e}_{cd} \cdot\tilde{\omega}_{\alpha\beta}(\tilde{\bf x})>  \\
& = <(-1)^{\gamma+\tau}\sum_{j=1}^n(x_{j\gamma}\tilde{\bf e}_{j\tau}-x_{j\tau}\tilde{\bf e}_{j\gamma}),\tilde{\bf e}_{ab}\otimes \tilde{\bf e}_{cd} \cdot \sum_{i=1}^n(-1)^{\alpha+\beta}(x_{i\alpha}\tilde{\bf e}_{i\beta}-x_{i\beta}\tilde{\bf e}_{i\alpha}) > \\
& = (-1)^{\alpha+\beta+\gamma+\tau}\sum_{i,j=1}^n(x_{j\gamma}x_{i\alpha}\delta_{ic}\delta_{d\beta}\delta_{ja}\delta_{b\tau}-
x_{j\gamma}x_{i\beta}\delta_{ic}\delta_{d\alpha}\delta_{ja}\delta_{b\tau} \\ 
& \quad -x_{j\tau}x_{i\alpha}\delta_{ic}\delta_{d\beta}\delta_{ja}\delta_{b\gamma}+
x_{j\tau}x_{i\beta}\delta_{ic}\delta_{d\alpha}\delta_{ja}\delta_{b\gamma})\\
& =  (-1)^{\alpha+\beta+\gamma+\tau}(x_{a\gamma}x_{c\alpha}\delta_{d\beta}\delta_{b\tau}-x_{a\gamma}x_{c\beta}\delta_{d\alpha}\delta_{b\tau}-x_{a\tau}x_{c\alpha}\delta_{d\beta}\delta_{b\gamma}+
x_{a\tau}x_{c\beta}\delta_{d\alpha}\delta_{b\gamma}) \\
&=  (-1)^{\alpha+\beta+\gamma+\tau}(x_{c\beta}\delta_{d\alpha}-x_{c\alpha}\delta_{d\beta})(x_{a\tau}\delta_{b\gamma}-x_{a\gamma}\delta_{b\tau}),
\end{align*}
where $\tilde{\omega}_{\alpha\beta}(\tilde{\bf x})=\mathcal{I}(X\Omega_{\alpha\beta})= (-1)^{\alpha+\beta}\sum_{i=1}^n(x_{i\alpha}\tilde{\bf e}_{i\beta}-x_{i\beta}\tilde{\bf e}_{i\alpha}),\,\,1\leq \alpha<\beta\leq n.$

Consequently, 
$$<\tilde{\omega}_{\gamma\tau}(\tilde{\bf x}), \tilde{\bf e}_{sj}\otimes \tilde{\bf e}_{sj} \cdot\tilde{\omega}_{\alpha\beta}(\tilde{\bf x})>= (-1)^{\alpha+\beta+\gamma+\tau}(x_{s\beta}\delta_{j\alpha}-x_{s\alpha}\delta_{j\beta})(x_{s\tau}\delta_{j\gamma}-x_{s\gamma}\delta_{j\tau}),$$
$$<\tilde{\omega}_{\gamma\tau}(\tilde{\bf x}), \tilde{\bf e}_{pj}\otimes \tilde{\bf e}_{qj} \cdot\tilde{\omega}_{\alpha\beta}(\tilde{\bf x})>= (-1)^{\alpha+\beta+\gamma+\tau}(x_{q\beta}\delta_{j\alpha}-x_{q\alpha}\delta_{j\beta})(x_{p\tau}\delta_{j\gamma}-x_{p\gamma}\delta_{j\tau}),$$
$$<\tilde{\omega}_{\gamma\tau}(\tilde{\bf x}), \tilde{\bf e}_{qj}\otimes \tilde{\bf e}_{pj} \cdot\tilde{\omega}_{\alpha\beta}(\tilde{\bf x})>= (-1)^{\alpha+\beta+\gamma+\tau}(x_{p\beta}\delta_{j\alpha}-x_{p\alpha}\delta_{j\beta})(x_{q\tau}\delta_{j\gamma}-x_{q\gamma}\delta_{j\tau}).$$

{\bf III.} The unit quaternions $\mathbf{q}=(q^0,q^1,q^2,q^3)\in S^3\subset \mathbb{R}^{4}$ and $-\mathbf{q}\in S^3\subset\mathbb{R}^{4}$ correspond to the following rotation in $SO(3)$:
\begin{equation}\label{Rq}
\mathbf{R}^{\mathbf{q}}=\left(
\begin{array}{ccc}
(q^0)^{2}+(q^{1})^{2}-(q^{2})^{2}-(q^{3})^{2} & 2(q^{1}q^{2}-q^{0}q^{3}) & 2(q^{1}q^{3}+q^{0}q^{2}) \\
2(q^{1}q^{2}+q^{0}q^{3})  & (q^0)^{2}-(q^{1})^{2}+(q^{2})^{2}-(q^{3})^{2}& 2(q^{2}q^{3}-q^{0}q^{1}) \\
2(q^{1}q^{3}-q^{0}q^{2}) & 2(q^{2}q^{3}+q^{0}q^{1}) & (q^0)^{2}-(q^{1})^{2}-(q^{2})^{2}+(q^{3})^{2}
\end{array}
\right).
\end{equation}

\bigskip

{\bf Acknowledgments.} This work was supported by the grant of the Romanian National Authority for
Scientific Research, CNCS - UEFISCDI, under the Romania-Cyprus bilateral cooperation programme
(module III), project number 760/2014. We are also thankful to Ioan Casu for his help with Maple programming.


\begin{thebibliography}{99}

\bibitem{absil-mahony-sepulchre-1} {\bf P.A. Absil, R. Mahony, R. Sepulchre}, {\it Optimization Algorithms on Matrix Manifolds}, Princeton University Press, 2008.
\bibitem{absil-mahony-trumpf} {\bf P.A. Absil, R. Mahony, J. Trumpf}, {\it An Extrinsic Look at the Riemannian Hessian}, Geometric Science of Information, Lecture Notes in Computer Science, Volume 8085, (2013), pp 361-368.
\bibitem{beck-hall}{\bf J.A. Beck, C.D. Hall}, {\it Relative equilibria of a rigid satellite in a circular
Keplerian orbit}, J. Astronaut. Sci., Vol. 40, Issue 3 (1998), pp. 215-247.
\bibitem{birtea-comanescu} {\bf P. Birtea, D. Com\u anescu}, {\it Geometric Dissipation for dynamical systems},
    Comm. Math. Phys., Vol. 316, Issue 2 (2012), pp. 375-394.
\bibitem{birtea-comanescu-popa} {\bf P. Birtea, D. Com\u anescu, C.A. Popa }, {\it Averaging on Manifolds by Embedding Algorithm},
    J. Math. Imaging Vis., Vol. 49, Issue 2 (2014), pp. 454-466.
\bibitem{comanescu} {\bf D. Com\u{a}nescu}, {\it The stability problem for the torque-free gyrostat investigated by using algebraic methods}, Applied Mathematics Letters, Volume 25, Issue 9 (2012), pp. 1185-1190.
\bibitem{comanescu-1} {\bf D. Com\u{a}nescu}, {\it Stability of equilibrium states in the Zhukovski case of heavy
gyrostat using algebraic methods}, Mathematical Methods in the Applied Sciences, Volume 36, Issue 4 (2013), pp. 373-382.
\bibitem{comanescu-2} {\bf D. Com\u{a}nescu}, {\it A note on stability of the vertical uniform rotations of the heavy top}, ZAMM, Volume 93, Issue 9 (2013), pp. 697-699.
\bibitem{donoho}{\bf D.L. Donoho, C. Grimes}, {\it Hessian eigenmaps: Locally linear embedding techniques for high-dimensional data}, Proceedings of the National Academy of Sciences, 100(10) (2003), pp. 5591-5596.
\bibitem{edelman} {\bf A. Edelman, T.A. Arias, S.T. Smith}, {\it The geometry of algorithms with orthogonality
    constraints}, SIAM J. Matrix Anal. Appl., Vol. 20, Issue 2 (1998), pp. 303-353.
\bibitem{ferreira}{\bf R. Ferreira, J. Xavier, J. P. Costeira, V. Barroso}, {\it Newton Algorithms for Riemannian Distance Related Problems on Connected Locally Symmetric Manifolds}, IEEE Journal of Selected Topics in Signal Processing, Volume 7, Issue 4 (2013), pp. 634-645.     
\bibitem{gallot}{\bf S. Gallot, D. Hulin, J. Lafontaine}, {\it Riemannian Geometry}, Universitext, Springer-Verlag, Berlin, 3rd edition, 2004.  
\bibitem{maddocks-1991}{\bf J.H. Maddocks}, {\it  Stability of relative equilibria}, IMA J. Appl. Math., Vol. 46 (1991), pp. 71-99.
\bibitem{moakher} {\bf M. Moakher}, {\it Means and averaging in the group of rotations}, SIAM J. Matrix Anal. Appl., Vol. 24, Issue 1 (2002), pp. 1-16.
\bibitem{montaldi} {\bf J. A. Montaldi, M. Rodriguez-Olmos}, {\it On the stability of Hamiltonian relative equilibria with non-trivial isotropy}, Nonlinearity, Vol. 24, Issue 10 (2011), pp. 2777-2783.
\bibitem{ortega-1}{\bf J.-P. Ortega, T. S. Ratiu}, {\it Stability of Hamiltonian relative equilibria}, Nonlinearity, Volume 12, Issue 3 (1999), pp. 693-720.
\bibitem{ortega-2}{\bf J.-P. Ortega, T. S. Ratiu}, {\it Non-linear stability of singular relative periodic orbits in Hamiltonian systems with symmetry}, Journal of Geometry and Physics, Volume 32, Issue 2 (1999), pp. 160-188.
\bibitem{patrick}{\bf G.W. Patrick}, {\it  Relative equilibria in Hamiltonian systems: the dynamic interpretation of nonlinear stability on a reduced phase space}, Journal of Geometry and Physics, Volume 9 (1992), pp. 111-119.
\bibitem{patrick-roberts}{\bf G.W. Patrick, M. Roberts, C. Wulff}, {\it  Stability of Poisson equilibria and Hamiltonian relative equilibria by energy methods}, Archive for Rational Mechanics and Analysis, Volume 174, Issue 3 (2004), pp. 301-344.
\bibitem{olmos}{\bf M. Rodriguez-Olmos}, {\it Stability of relative equilibria with singular momentum values in simple mechanical systems}, Nonlinearity, Volume 19, Issue 4 (2006).
\bibitem{wang-xu}{\bf Y. Wang, S. Xu}, {\it Equilibrium attitude and nonlinear attitude stability of a spacecraft on a stationary orbit
around an asteroid}, J. Adv. Space Res., Vol. 52, Issue 8 (2013), pp. 1497-1510.

\end{thebibliography}
\end{document}